# Automated tuning of bifurcations via feedback

Luc Moreau[*]     Eduardo Sontag[†]     Murat Arcak[‡]

October 24, 2018

## 1 Introduction

The present paper studies a feedback regulation problem that arises in at least two different biological applications. The feedback regulation problem under consideration may be interpreted as an adaptive control problem, but has not yet been studied in the control literature. The goal of the paper is to introduce this problem and to present some first results.

The feedback regulation problem is concerned with a forced dynamical system

$$\dot{x} = f_\mu(x, u(t)) \tag{1}$$

which depends on a parameter $\mu$. The input $u(t)$ in equation (1) represents an external stimulus; it is not a control variable. We are interested in the input-to-state properties of equation (1). In general, the input-to-state properties of equation (1) depend on the parameter $\mu$. In particular, the input-to-state properties may change drastically as a function of $\mu$ when the unforced dynamics exhibit a bifurcation. In this case, quite interesting and remarkable amplification properties may result from operating in the proximity of the bifurcation point. We illustrate this with two biological examples from the literature.

The first example is concerned with the auditory system. In order to detect the sounds of the outside world, hair cells in the cochlea operate as nanosensors which transform acoustic stimuli into electric signals. In a series of recent

[*]Postdoctoral Fellow of the Fund for Scientific Research - Flanders (Belgium) (F.W.O.-Vlaanderen) and recipient of an Honorary Fellowship of the Belgian American Educational Foundation. This paper presents research results of the Belgian Programme on Inter-University Poles of Attraction, initiated by the Belgian State, Prime Minister's Office for Science, Technology and Culture. The scientific responsibility rests with its authors. SYS-TeMS, Ghent University, Technologiepark 9, 9052 Zwijnaarde, Belgium. Currently visiting Mechanical and Aerospace Engineering Department, Princeton University, Princeton N.J. 08544, USA

[†]Department of Mathematics, Rutgers, The State University of New Jersey, Piscataway N.J. 08854, USA. This work was supported in part by US Air Force Grant F49620-01-1-0063, and by National Institutes of Health Grant P20 GM64375.

[‡]Dept. of Electrical, Computer and Systems Engineering, Rensselaer Polytechnic Institute, Troy, NY 12180, USA

papers [4, 5, 6, 7], the hair cells in the cochlea are modeled as active, almost self-oscillating systems. Ions such as $Ca^{++}$ are believed to contribute to the hair cell's tendency to exhibit spontaneous oscillations. For low concentrations of the ions, the viscous damping forces of the fluid that surrounds the hair cells dominate and the hair cell oscillations are damped. As the concentration increases the system undergoes a Hopf bifurcation, the dynamics become unstable and the hair cells exhibit spontaneous oscillations. It is argued in these papers that the hair cells operate in the proximity of this bifurcation point, where the activity of the ions compensates the damping effects. In this case, even a weak sound stimulus can cause a detectable oscillation. (This follows from the generic properties of a forced dynamical system exhibiting a Hopf bifurcation). This mechanism thus provides an explanation for the experimentally observed ultrahigh sensitivity of the ear.

The second example arises in the study of persistent neural activity [3, 12, 14, 13]. Neural activity of a single neuron has a natural tendency to decay with a relaxation time of about 5–100 ms. This natural tendency to decay can be opposed by positive synaptic feedback loops. If this feedback is weak, then the natural tendency to decay dominates and neural activity still decreases. As the feedback gain is increased, the neural dynamics undergo a bifurcation and the dynamics become unstable. When the feedback is tuned to exactly balance the decay, then neural activity neither increases nor decreases but persists without change. According to a long-standing hypothesis, this is the mechanism that so-called neural integrators use to maintain persistent neural activity. A transient stimulus of a neural integrator can then cause a persistent change in neural activity. This mechanism forms the basis for short-term analogue memory and plays a central role in the oculomotor control system.

Both examples illustrate the remarkable and interesting input-to-state properties of a dynamical system poised at a bifurcation. It is clear, however, that operating in the proximity of a bifurcation point requires a fine-tuning of parameters. And thus the question arises as to how a (biological) system can be tuned with high precision to its critical bifurcation point. In both examples it has been suggested that feedback regulation of the bifurcation parameter may provide a robust mechanism to ensure the required fine-tuning of parameters [4, 6, 3, 16]. In the literature on hearing this has given rise to the terminology of "self-tuned Hopf bifurcation". Quite remarkably both studies have been pursued independently of each other in spite of their strong similarities which are obvious from a control engineering perspective.

The present paper aims to initiate a mathematical study of this feedback regulation problem. The problem may be formulated as follows. Find an adaptation law which steers the bifurcation parameter $\mu$ to its critical value $\mu_0$ without prior knowledge of this critical value $\mu_0$. Our main motivation for studying this problem stems from the fact that the precise biophysical nature of the feedback mechanisms involved in the above biological applications is still unknown. A theoretical study of this problem may guide the search for possible biophysical mechanisms. Even when the precise biophysical adaptation mechanisms will have been discovered, it is to be expected that a profound un-derstanding of the mechanisms involved can only be achieved when insight has been gained in the fundamental trade-offs and limitations that may be inherent to this problem. This is very similar in spirit to the internal model principle which has turned out to be central to the understanding of adaptation in bacterial chemotaxis [17], and it may be seen as an illustration of the important role that control engineering may play in the emerging field of systems biology.

A second motivation for the present study is of a more fundamental nature. In view of the quite remarkable and interesting input-to-state properties that arise when operating in the proximity of a bifurcation point, the problem under consideration may turn out to be of independent interest. Despite its mathematical appeal, it has not attracted attention before in the control community. The present paper brings together several well-established techniques from nonlinear and adaptive control to derive some first results for this problem.

We end this introduction with some references to related work. The present research bears some similarity with the problem of experimental instability detection [2], where an operating parameter is adapted online in order to efficiently locate bifurcations through experiments. A different problem which is related to the present study concerns the detection and prediction of instabilities via closed-loop monitoring techniques [9, 10]. The possible implications of the present work for the problems studied in these papers remains a topic for further research.

## 2  Self-tuning of a first-order system

The one-dimensional system

$$\dot{x} = (\mu - \mu_0)x + u(t) \tag{2}$$

captures some of the essential features of the neural integrator. In this interpretation $\mu_0$ represents the natural decay rate of neural activity and $\mu$ corresponds to the synaptic feedback gain. We view $\mu$ as an adjustable parameter and $\mu_0$ as an unknown constant. Clearly this system exhibits a bifurcation. If $\mu = \mu_0$ then equation (2) behaves as a perfect integrator, if $\mu < \mu_0$ (respectively $\mu > \mu_0$) then equation (2) is referred to as a *leaky* (respectively *unstable*) integrator. We ask the question as to how proximity to the bifurcation point may be ensured.

We study this question from an engineering perspective. In a first approximation we ignore the presence of the input and ask the following simpler question. Find an adaptation law for the parameter $\mu$ which steers $\mu$ to its bifurcation value $\mu_0$ for the system

$$\dot{x} = (\mu - \mu_0)x. \tag{3}$$

This adaptation law may depend on $x$ and $\mu$, but should be independent of $\mu_0$, as this value is not known (precisely). Let us first discuss the feasibility of this problem. It is easy to see that if $x = 0$ at some time instant then $x = 0$ for all times and in this case it is clearly impossible to steer $\mu$ to $\mu_0$ without prior knowledge of $\mu_0$. We therefore restrict attention to the set of strictly positive values for $x$, which is invariant under the dynamics (3). Considering only strictly positive values for $x$ is physically relevant, as this variable represents a level of neural activity (rate of action potential firing).

The following theorem provides sufficient conditions for the adaptation law

$$\dot{\mu} = f(x) - g(\mu) \tag{4}$$

to steer $\mu$ to its bifurcation value $\mu_0$ for the system (3).

**Theorem 1.** *Let $\mu_0 \in \mathbb{R}$ and consider continuously differentiable functions $f : \mathbb{R}_{>0} \to \mathbb{R}$ and $g : \mathbb{R} \to \mathbb{R}$. Assume that $f$ is strictly decreasing, $g$ is strictly increasing, and $g(\mu_0)$ is in the image of $f$. Then the nonlinear system (3)–(4) with $x \in \mathbb{R}_{>0}$ and $\mu \in \mathbb{R}$ has a unique equilibrium point, which is globally asymptotically stable (and locally exponentially stable if $\mathrm{d}f/\mathrm{d}x$ takes only strictly negative values and $\mathrm{d}g/\mathrm{d}\mu$ only strictly positive values). At equilibrium $\mu$ equals $\mu_0$.*

The proof relies on a coordinate transformation which converts (3)–(4) into a nonlinear mass-spring-damper system. Global asymptotic stability follows readily from LaSalle's theorem. Local exponential stability is shown by means of the linearization principle.

*Proof.* We introduce new coordinates $q = \ln(x) - \ln(f^{-1}(g(\mu_0)))$ and $p = \mu - \mu_0$. This transformation from $(x, \mu)$ to $(q, p)$ is a global $C^\infty$-diffeomorphism from $\mathbb{R}_{>0} \times \mathbb{R}$ to $\mathbb{R}^2$. Expressed in these new coordinates (3)–(4) becomes

$$\dot{q} = p, \tag{5}$$
$$\dot{p} = f(\exp(q) f^{-1}(g(\mu_0))) - g(p + \mu_0). \tag{6}$$

The system of equations (5)–(6) has a unique equilibrium point at the origin. This equilibrium corresponds to an equilibrium in original coordinates where $\mu = \mu_0$.

First we prove that the null-solution of (5)–(6) is globally asymptotically stable. We rewrite (5)–(6) as

$$\dot{q} = p, \tag{7}$$
$$\dot{p} = -\tilde{f}(q) - \tilde{g}(p), \tag{8}$$

with $\tilde{f}(q) = -f(\exp(q) f^{-1}(g(\mu_0))) + g(\mu_0)$ and $\tilde{g}(p) = g(p + \mu_0) - g(\mu_0)$. Clearly $\tilde{f}$ and $\tilde{g}$ are both strictly increasing, continuously differentiable functions from $\mathbb{R}$ to $\mathbb{R}$ which are zero at zero. The candidate Lyapunov function

$$V : \mathbb{R}^2 \to \mathbb{R} : (q, p) \mapsto \int_0^q \tilde{f}(\xi) \, \mathrm{d}\xi + \frac{p^2}{2} \tag{9}$$

is twice continuously differentiable, positive definite and radially unbounded. Its time derivative along the solutions of (7)–(8) satisfies

$$\dot{V}(q,p) = -p\tilde{g}(p) \begin{cases} < 0 & \text{if } p \neq 0, \\ = 0 & \text{if } p = 0. \end{cases} \qquad (10)$$

Since the null-solution is the only solution of (7)–(8) along which $p$ vanishes identically, global asymptotic stability of the origin follows from LaSalle's theorem [8, Corollary 3.2].

Next we prove that the equilibrium of (3)–(4) is locally exponentially stable if $df/dx$ takes only strictly negative values and $dg/d\mu$ only strictly positive values. Since the transformation from $(x, \mu)$ to $(q, p)$ is a $C^\infty$-diffeomorphism it suffices to prove that the null-solution of (5)–(6) is locally exponentially stable. This follows readily from the linearization principle. The linearization of (5)–(6) around the origin is given by

$$\dot{q} = p, \qquad (11)$$
$$\dot{p} = \frac{df}{dx}(f^{-1}(g(\mu_0)))f^{-1}(g(\mu_0))q - \frac{dg}{d\mu}(\mu_0)p. \qquad (12)$$

If $df/dx$ takes only strictly negative values and $dg/d\mu$ only strictly positive values then the linearization (11)–(12) is exponentially stable, as required. □

## Discussion of Theorem 1

1. A mathematically appealing adaptation law which satisfies the assumptions of Theorem 1 for all possible values of $\mu_0$ is given by

$$\dot{\mu} = -a\ln(x) - b\mu \qquad (13)$$

with $a, b > 0$. This adaptation law ensures convergence to the bifurcation point for any possible value of $\mu_0$. It has the interesting property that the coordinate transformation introduced in the proof of Theorem 1 converts (3), (13) into the *linear* system

$$\dot{q} = p, \qquad (14)$$
$$\dot{p} = -aq - bp. \qquad (15)$$

From a biological point of view however this adaptation law is probably not very relevant, among others because the right hand side of (13) is not bounded. A bounded adaptation law which could be more relevant for biological applications is for example

$$\dot{\mu} = 1/(1+x^2) - 1/(1+\exp(-\mu)). \qquad (16)$$

This bounded adaptation law also ensures convergence to the bifurcation point for any possible value of $\mu_0$.

2. There is a subtle relationship between self-tuning of bifurcations and ideas from robust control theory. This relation is made explicit by equation (5), which represents an integrator and corresponds to integral action studied in robust control theory. Although perhaps surprising at first sight, this relation may be understood by regarding the constant $\mu_0$ as an unknown perturbation acting on the system. It is well-known from robust control theory that (under appropriate conditions) rejecting or tracking an unknown, constant disturbance requires integral action. (This is a special case of the internal model principle—see [17, 15] for a discussion of the internal model principle from a systems biology perspective). The present manifestation of integral action, however, differs from traditional robust control situations because the integral action is not generated by the adaptation law but is inherent to the dynamics of system (3) itself.

3. If an additional perturbation term $\varepsilon p(x, \mu, t)$ is considered

$$\dot{x} = (\mu - \mu_0)x + \varepsilon p(x, \mu, t), \qquad (17)$$

we may conclude that the adaptation law (4) steers $\mu$ approximately to $\mu_0$ under the conditions of Theorem 1 provided the perturbation is small enough.

**Theorem 2.** *Let $\mu_0 \in \mathbb{R}$ and consider continuously differentiable functions $f : \mathbb{R}_{>0} \to \mathbb{R}$ and $g : \mathbb{R} \to \mathbb{R}$. Assume that $f$ is strictly decreasing, $g$ is strictly increasing, and $g(\mu_0)$ is in the image of $f$. Consider a strictly positive parameter $\varepsilon$ and a continuous function $p : \mathbb{R}_{>0} \times \mathbb{R}^2 \to \mathbb{R} : (x, \mu, t) \to p(x, \mu, t)$ and assume that $p$ is bounded in $t$ uniformly with respect to $(x, \mu)$ belonging to compact subsets of $\mathbb{R}_{>0} \times \mathbb{R}$. Then the unique equilibrium point for the original system (3)–(4) is uniformly semiglobally practically asymptotically stable for the perturbed system (17)–(4).*

(We refer to the Appendix for a definition of uniform semiglobal practical asymptotic stability.) In particular, this theorem guarantees that along all trajectories of (17)–(4) originating in an arbitrarily large compact subset of $\mathbb{R}_{>0} \times \mathbb{R}$, $\mu$ will eventually converge to an arbitrarily small open neighborhood of $\mu_0$ provided $\varepsilon$ is small enough. This result is a manifestation of the well-known robustness of asymptotic stability with respect to small perturbations. It provides some justification for ignoring the external stimulus in equation (2) and studying the simpler equation (3) instead, at least in a first approximation.

## 3 A second-order system

In a series of articles [4, 5, 6, 7] it is argued that various nonlinear phenomena in the auditory system (such as ultrahigh sensitivity to weak signals) may be understood in terms of the generic properties of a forced dynamical system exhibiting a Hopf bifurcation. In [6] this is illustrated by means of a standard model for nonlinear oscillations which (with the notation of the present paper) takes the form

$$\ddot{x} + (\mu_0 - \mu)\dot{x} + \lambda\dot{x}^3 + \omega^2 x = u(t). \tag{18}$$

As before $\mu$ is an adjustable parameter and $\mu_0$ is an unknown constant. In the absence of an external input $u$, equation (18) describes damped oscillations if $\mu < \mu_0$. At $\mu = \mu_0$ the system undergoes a Hopf bifurcation and for $\mu > \mu_0$ spontaneous oscillations are being generated. Assuming that it is possible to tune $\mu = \mu_0$ such a model captures the biophysical nature of hair cell oscillations within the *cochlea* where the hair cells are thought of as active, almost self-oscillating systems rather than passive oscillators. In [4, 5, 6, 7] it is shown how such a dynamical system operating at a Hopf bifurcation indeed may account for several of the observed nonlinear phenomena in hearing. Again the question arises as to how proximity to the bifurcation point may be ensured. In [4, 6] it is suggested that a feedback mechanism is responsible for this self-tuning and numerical simulations are provided to support this claim.

Here we want to contribute to a mathematical analysis of self-tuning of an oscillator. In a first approximation we ignore the cubic damping term and the external input and we study the following question. Find an adaptation law for the parameter $\mu$ which does not depend on $\mu_0$ and which steers $\mu$ to its bifurcation value $\mu_0$ for the system

$$\ddot{x} + (\mu_0 - \mu)\dot{x} + \omega^2 x = 0. \tag{19}$$

As before let us first discuss the feasibility of this problem. It is easy to see that if $x = \dot{x} = 0$ at some time instant then $x = \dot{x} = 0$ for all times and in this case it is clearly impossible to steer $\mu$ to $\mu_0$ without prior knowledge of $\mu_0$. We therefore restrict attention to the set where $x$ and $\dot{x}$ are not both zero—this set is invariant under the dynamics (19).

The following theorem asserts that the adaptation law

$$\dot{\mu} = -a\ln(\sqrt{x^2 + (\dot{x}/\omega)^2}) - b\mu \tag{20}$$

steers $\mu$ to its bifurcation value $\mu_0$ for the system (19) if $0 < a \leq b^2$ and $b > 0$.

**Theorem 3.** *Let $\mu_0 \in \mathbb{R}$ and $a, b, \omega \in \mathbb{R}_{>0}$. Assume that $a \leq b^2$. The nonlinear system (19)–(20) with $(x, \dot{x}) \in \mathbb{R}^2 \setminus \{(0,0)\}$ and $\mu \in \mathbb{R}$ has a unique periodic orbit[1] which is globally asymptotically stable and locally exponentially stable. On this periodic orbit $\mu = \mu_0$.*

The proof of global asymptotic stability is based on passivity techniques. We interpret (19)–(20) (in a different coordinate system) as a feedback interconnection of a linear system with a dynamic feedback which satisfies a sector condition. Via the Kalman-Yakubovich-Popov lemma we obtain a proper, non-increasing Lyapunov function. Global asymptotic stability of the periodic orbit then follows from LaSalle's invariance principle.

---

[1] A periodic orbit is a subset of the state-space which is the image of a periodic solution.

*Proof.* First we introduce new coordinates $r$ and $\phi$ according to $x = r\cos(\phi)$ and $\dot{x} = -r\omega\sin(\phi)$. The transformation from $(x, \dot{x})$ to $(r, \phi)$ is a global $C^\infty$-diffeomorphism from $\mathbb{R}^2 \setminus \{(0,0)\}$ to $\mathbb{R}_{>0} \times S^1$. Expressed in the coordinates $(r, \phi, \mu)$, equations (19)–(20) become

$$\dot{r} = (\mu - \mu_0) r \sin^2(\phi), \tag{21}$$
$$\dot{\phi} = \omega + (\mu - \mu_0) \sin(\phi) \cos(\phi), \tag{22}$$
$$\dot{\mu} = -a \ln(r) - b\mu. \tag{23}$$

Next we introduce new coordinates $q = \ln(r) + b\mu_0/a$ and $p = \mu - \mu_0$. The transformation from $(r, \mu)$ to $(q, p)$ is a global $C^\infty$-diffeomorphism from $\mathbb{R}_{>0} \times \mathbb{R}$ to $\mathbb{R}^2$. Expressed in the coordinates $(q, \phi, p)$, equations (21)–(23) become

$$\dot{q} = p \sin^2(\phi), \tag{24}$$
$$\dot{\phi} = \omega + p \sin(\phi) \cos(\phi), \tag{25}$$
$$\dot{p} = -aq - bp. \tag{26}$$

Clearly the system of equations (24)–(26) has a periodic orbit $\{(q, \phi, p) \in \mathbb{R} \times S^1 \times \mathbb{R} : q = p = 0\}$ which we denote by $\mathcal{A}$. This periodic orbit corresponds to a periodic orbit in original coordinates where $\mu = \mu_0$. It is clear from the following paragraph that this periodic orbit is unique.

First we prove that the periodic orbit $\mathcal{A}$ of (24)–(26) is globally asymptotically stable. We interpret (24)–(26) as a feedback interconnection of a linear control system

$$\dot{q} = -u, \tag{27}$$
$$\dot{p} = -aq - bp, \tag{28}$$

with a negative, dynamic feedback $u = -\sin^2(\phi) p$ where $\phi$ satisfies (25). The transfer function $H(s)$ of (27)–(28) from $u$ to $p$ is given by

$$H(s) = \frac{a}{s(s+b)}. \tag{29}$$

Since $a$ and $b$ are strictly positive and $a \leq b^2$ the transfer function $H(s) + 1$ is positive real. Hence by (a modified version of) the Kalman-Yakubovich-Popov lemma (proven in [1] and described in [8, Exercise 10.2]) there exists a positive definite, quadratic Lyapunov function $(q, p) \mapsto V(q, p)$ whose time derivative along the solutions of the control system (27)–(28) satisfies

$$\dot{V}(q, p, u) \leq up + u^2. \tag{30}$$

In other words the control system (27)–(28) with input $u$ and output $p$ is input feedforward passive with a shortage of passivity [11]. Since the feedback satisfies

$$up + u^2 = (\sin^4(\phi) - \sin^2(\phi)) p^2 = -(\sin(\phi)\cos(\phi) p)^2, \tag{31}$$

it follows that the time-derivative of $V$ along the solutions of (24)–(26) satisfies

$$\dot{V}(q,\phi,p) \begin{cases} < 0 & \text{if } \sin(\phi)\cos(\phi)p \neq 0, \\ = 0 & \text{if } \sin(\phi)\cos(\phi)p = 0. \end{cases} \qquad (32)$$

Since $(q,\phi,p) \mapsto V(q,p)$ is positive definite with respect to the periodic orbit $\mathcal{A}$ it follows that $\mathcal{A}$ is stable. In addition, since $(q,\phi,p) \mapsto V(q,p)$ is radially unbounded with respect to $\mathcal{A}$ and since $\mathcal{A}$ is compact it follows that every solution of (24)–(26) is bounded. Finally, since $\mathcal{A}$ is the largest invariant set of (24)–(26) contained in $\{(q,\phi,p) \in \mathbb{R} \times S^1 \times \mathbb{R} : \dot{V}(q,\phi,p) = 0\}$ it follows from LaSalle's theorem that every solution of (24)–(26) converges to $\mathcal{A}$.

It remains to be proven that the periodic orbit of (19)–(20) is locally exponentially stable. This follows directly from Theorem 4 which is stated and proven below. $\square$

Although mathematically appealing, the adaptation law (20) is probably not very relevant from a biological point of view, among others because the right hand side is not bounded. It would therefore be interesting to have a result available that applies to more general adaptation laws

$$\dot{\mu} = f(\sqrt{x^2 + (\dot{x}/\omega)^2}) - g(\mu). \qquad (33)$$

This is the subject of Theorem 4. Unlike the previous results, Theorem 4 is a local stability result.

**Theorem 4.** *Let $\mu_0 \in \mathbb{R}$ and consider continuously differentiable functions $f : \mathbb{R}_{>0} \to \mathbb{R}$ and $g : \mathbb{R} \to \mathbb{R}$. Assume that $g(\mu_0)$ is in the image of $f$. Consider $r^* \in f^{-1}(g(\mu_0))$ and assume that $0 < -(\mathrm{d}f/\mathrm{d}r)(r^*)r^* \leq ((\mathrm{d}g/\mathrm{d}\mu)(\mu_0))^2$ and $(\mathrm{d}g/\mathrm{d}\mu)(\mu_0) > 0$. Then the system of equations (19) and (33) with $(x,\dot{x}) \in \mathbb{R}^2 \setminus \{(0,0)\}$ and $\mu \in \mathbb{R}$ has a periodic orbit which is locally exponentially stable and where $\mu = \mu_0$.*

*Proof.* We proceed along the lines of the proof of Theorem 3. First we introduce new coordinates $r$ and $\phi$ according to $x = r\cos(\phi)$ and $\dot{x} = -r\omega\sin(\phi)$. The transformation from $(x,\dot{x})$ to $(r,\phi)$ is a global $C^\infty$-diffeomorphism from $\mathbb{R}^2 \setminus \{(0,0)\}$ to $\mathbb{R}_{>0} \times S^1$. Expressed in the coordinates $(r,\phi,\mu)$, equations (19) and (33) become

$$\dot{r} = (\mu - \mu_0)r\sin^2(\phi), \qquad (34)$$
$$\dot{\phi} = \omega + (\mu - \mu_0)\sin(\phi)\cos(\phi), \qquad (35)$$
$$\dot{\mu} = f(r) - g(\mu). \qquad (36)$$

Next we introduce new coordinates $q = \ln(r) - \ln(r^*)$ and $p = \mu - \mu_0$. The transformation from $(r,\mu)$ to $(q,p)$ is a global $C^\infty$-diffeomorphism from $\mathbb{R}_{>0} \times \mathbb{R}$ to $\mathbb{R}^2$. Expressed in the coordinates $(q,\phi,p)$, equations (34)–(36) become

$$\dot{q} = p\sin^2(\phi), \qquad (37)$$
$$\dot{\phi} = \omega + p\sin(\phi)\cos(\phi), \qquad (38)$$
$$\dot{p} = f(\exp(q)r^*) - g(p + \mu_0). \qquad (39)$$

Clearly the system of equations (37)–(39) has a periodic orbit $\{(q, \phi, p) \in \mathbb{R} \times S^1 \times \mathbb{R} : q = p = 0\}$ which we denote by $\mathcal{A}$. This periodic orbit corresponds to a periodic orbit in original coordinates where $\mu = \mu_0$.

Since the transformation from $(x, \dot{x}, \mu)$ to $(q, \phi, p)$ is a $C^\infty$-diffeomorphism it suffices to prove that the periodic orbit $\mathcal{A}$ of (37)–(39) is locally exponentially stable. We first ignore the "higher order" terms in the right hand side of (38) and (39) and consider the simpler system

$$\dot{q} = p \sin^2(\phi), \tag{40}$$
$$\dot{\phi} = \omega, \tag{41}$$
$$\dot{p} = -aq - bp, \tag{42}$$

where we have introduced the notation $a = -(df/dr)(r^*)r^*$ and $b = (dg/d\mu)(\mu_0)$. By the assumptions of the theorem $0 < a \leq b^2$ and $b > 0$. Repeating the arguments of the proof of Theorem 4 it is easy to see that the set $\mathcal{A}$ is a globally asymptotically stable periodic orbit of (40)–(42). Because of the special structure of equation (41) we may interpret $\phi$ as a time variable and the system of equations (40) and (42) as a periodically time-varying linear system. With this interpretation the periodic orbit $\mathcal{A}$ of (40)–(42) corresponds to the null-solution of the linear system of equations (40) and (42). Standard converse theorems for periodically time-varying linear systems (see for example [8, Theorem 3.10]) yield the existence of strictly positive real numbers $c_1$, $c_2$ and $c_3$ and a continuously differentiable Lyapunov function $(q, \phi, p) \mapsto V(q, \phi, p)$ which is quadratic in $(q, p)$ and (together with its partial derivative $\partial V/\partial \phi$) bounded in $\phi$ such that

$$c_1(q^2 + p^2) \leq V(q, \phi, p) \leq c_2(q^2 + p^2) \tag{43}$$

and such that the time derivative of $V$ along the solutions of (40)–(42) satisfies

$$\dot{V}(q, \phi, p) \leq -c_3(q^2 + p^2). \tag{44}$$

With this Lyapunov function $V$ we now prove that the periodic orbit $\mathcal{A}$ of (37)–(39) is locally exponentially stable. Indeed, for $q$ and $p$ sufficiently close to zero the time derivative of $V$ evaluated along the solutions of (37)–(39) satisfies

$$\dot{V}(q, \phi, p) \leq -\frac{c_3}{2}(q^2 + p^2) \tag{45}$$

since the extra terms in the right hand side of (38) and (39) give rise to extra terms in the Lyapunov balance which are bounded in $\phi$ and of order higher than two in $(q, p)$. This shows that the periodic orbit $\mathcal{A}$ of (37)–(39) is locally exponentially stable. □

## Discussion of Theorems 3 and 4

1. The assumptions of Theorem 4 involve the unknown critical value $\mu_0$. If $\mu_0$ is known to belong to some interval, it is of interest to have an adaptation law which satisfies the assumptions of the theorem for all $\mu_0$ in this interval. For example, it is easily verified that the adaptation law

$$\dot{\mu} = 1/(1 + \left(\sqrt{x^2 + (\dot{x}/\omega)^2}\right)^a) - b\mu \tag{46}$$

with $0 < a \leq 4b^2$ and $b > 0$ satisfies all assumptions of the theorem if $0 < \mu_0 < 1/b$.

2. The proposed adaptation laws (20) and (33) depend not only on $x$ and $\mu$ but also on $\dot{x}$ and $\omega$. Their implementation requires that the state variable $\dot{x}$ is measured and the parameter $\omega$ is known. This may be a disadvantage.

3. If an additional perturbation term $\varepsilon p(x, \dot{x}, \mu, t)$ is considered

$$\ddot{x} + (\mu_0 - \mu)\dot{x} + \omega^2 x = \varepsilon p(x, \dot{x}, \mu, t), \tag{47}$$

we may conclude that the adaptation law (20) or (33) steers $\mu$ approximately to $\mu_0$ under the conditions of Theorem 3 or Theorem 4 provided the perturbation is small enough.

**Theorem 5.** *Let $\mu_0 \in \mathbb{R}$ and consider continuously differentiable functions $f : \mathbb{R}_{>0} \to \mathbb{R}$ and $g : \mathbb{R} \to \mathbb{R}$. Assume that $g(\mu_0)$ is in the image of $f$. Consider $r^* \in f^{-1}(g(\mu_0))$ and assume that $0 < -(\mathrm{d}f/\mathrm{d}r)(r^*)r^* \leq ((\mathrm{d}g/\mathrm{d}\mu)(\mu_0))^2$ and $(\mathrm{d}g/\mathrm{d}\mu)(\mu_0) > 0$. Consider a strictly positive parameter $\varepsilon$ and a continuous function $p : (\mathbb{R}^2 \setminus \{0\}) \times \mathbb{R}^2 \to \mathbb{R} : (x, \dot{x}, \mu, t) \to p(x, \dot{x}, \mu, t)$ and assume that $p$ is bounded in $t$ uniformly with respect to $(x, \dot{x}, \mu)$ belonging to compact subsets of $(\mathbb{R}^2 \setminus \{0\}) \times \mathbb{R}$. Then there is a subset of $(\mathbb{R}^2 \setminus \{0\}) \times \mathbb{R}$ where $\mu = \mu_0$ and which is uniformly practically asymptotically stable for the perturbed system (47) and (33). If $f(\cdot)$ and $g(\mu)$ take the particular form $-a\ln(\cdot)$ respectively $b\mu$, then there is a subset of $(\mathbb{R}^2 \setminus \{0\}) \times \mathbb{R}$ where $\mu = \mu_0$ and which is uniformly semiglobally practically asymptotically stable.*

(We refer to the Appendix for a definition of the notion of uniform (semiglobal) practical asymptotic stability.) This result is a manifestation of the well-known robustness of asymptotic (or exponential) stability with respect to small perturbations. It provides some justification for ignoring the cubic nonlinearity and external stimulus in equation (18) and studying the simpler equation (19) instead, at least in a first approximation.

# A  Practical and semiglobal stability definitions

Let $\Omega$ be an open subset of $\mathbb{R}^n$ ($n \in \mathbb{N}$) and consider a family of continuous functions $f_\varepsilon : \mathbb{R} \times \Omega \to \mathbb{R}^n$, labeled by a parameter $\varepsilon > 0$. We are interested in the stability properties of the family of differential equations

$$\dot{x} = f_\varepsilon(t, x) \tag{48}$$

for small values of $\varepsilon$. Let $\mathcal{A}$ be a compact subset of $\Omega$, which need not necessarily be forward invariant for (48).

**Definition 1.** *For the family of differential equations* (48), *the set $\mathcal{A}$ is:*

1. *Uniformly practically stable if for every open neighborhood $U_2$ of $\mathcal{A}$ there is $\varepsilon^* > 0$ and an open neighborhood $U_1$ of $\mathcal{A}$ such that for all $\varepsilon \in (0, \varepsilon^*]$ every solution $\xi$ of* (48) *satisfies: if $\xi(t_0) \in U_1$ for some $t_0$ in the domain of $\xi$, then $\xi(t) \in U_2$ for all $t \geq t_0$ in the domain of $\xi$.*

2. *Uniformly practically asymptotically stable if it is uniformly practically stable and, in addition, there is an open neighborhood $U_1$ of $\mathcal{A}$ such that for every open neighborhood $U_2$ of $\mathcal{A}$ there is $T \geq 0$ and $\varepsilon^* > 0$ such that for all $\varepsilon \in (0, \varepsilon^*]$ every solution $\xi$ of* (48) *satisfies: if $\xi(t_0) \in U_1$ for some $t_0$ in the domain of $\xi$, then $\xi(t) \in U_2$ for all $t \geq t_0 + T$ in the domain of $\xi$.*

3. *Uniformly semiglobally bounded if for every compact subset $K_1$ of $\Omega$ there is $\varepsilon^* > 0$ and a compact subset $K_2$ of $\Omega$ such that for all $\varepsilon \in (0, \varepsilon^*]$ every solution $\xi$ of* (48) *satisfies: if $\xi(t_0) \in K_1$ for some $t_0$ in the domain of $\xi$, then $\xi(t) \in K_2$ for all $t \geq t_0$ in the domain of $\xi$.*

4. *Uniformly semiglobally practically asymptotically stable if it is uniformly practically stable and uniformly semiglobally bounded and, in addition, for every compact subset $K$ of $\Omega$ and every open neighborhood $U$ of $\mathcal{A}$ there is $T \geq 0$ and $\varepsilon^* > 0$ such that for all $\varepsilon \in (0, \varepsilon^*]$ every solution $\xi$ of* (48) *satisfies: if $\xi(t_0) \in K$ for some $t_0$ in the domain of $\xi$, then $\xi(t) \in U$ for all $t \geq t_0 + T$ in the domain of $\xi$.*